\documentclass[12pt,a4paper]{article}
\usepackage{textcomp}
\usepackage{listings}
\usepackage[space=true]{accsupp}
\newcommand{\copyablespace}{\BeginAccSupp{method=hex,unicode,ActualText=00A0}\EndAccSupp{}}
\usepackage[T1]{fontenc}
\usepackage[english]{babel}
\usepackage[colorlinks = true,
            linkcolor = blue,
            urlcolor  = blue,
            citecolor = blue,
            anchorcolor = blue]{hyperref}
\usepackage{graphicx}
\usepackage{amsmath}
\usepackage{amssymb}
\usepackage{mathtools}
\usepackage[table,dvipsnames]{xcolor}
\usepackage{fancyhdr}
\usepackage{longtable}
\usepackage{multicol}
\usepackage{enumerate}
\usepackage{enumitem}
\setlist[itemize]{leftmargin=5.5mm}
\usepackage{multirow}
\usepackage{natbib}
\usepackage{caption}
\usepackage{subcaption}
\usepackage{xurl}
\usepackage[utf8]{inputenc}
\usepackage{tikz}
\usepackage{pgfplots}
\usepgfplotslibrary{colorbrewer}
\pgfplotsset{compat=newest, cycle list/Set1-8}

\usetikzlibrary{pgfplots.statistics, pgfplots.colorbrewer}
\usepackage{pgfplotstable}

\pgfplotsset{
    box plot/.style={
        /pgfplots/.cd,
        black,
        only marks,
        mark=-,
        mark size=\pgfkeysvalueof{/pgfplots/box plot width},
        /pgfplots/error bars/y dir=plus,
        /pgfplots/error bars/y explicit,
        /pgfplots/table/x index=\pgfkeysvalueof{/pgfplots/box plot x index},
    },
    box plot box/.style={
        /pgfplots/error bars/draw error bar/.code 2 args={%
            \draw  ##1 -- ++(\pgfkeysvalueof{/pgfplots/box plot width},0pt) |- ##2 -- ++(-\pgfkeysvalueof{/pgfplots/box plot width},0pt) |- ##1 -- cycle;
        },
        /pgfplots/table/.cd,
        y index=\pgfkeysvalueof{/pgfplots/box plot box top index},
        y error expr={
            \thisrowno{\pgfkeysvalueof{/pgfplots/box plot box bottom index}}
            - \thisrowno{\pgfkeysvalueof{/pgfplots/box plot box top index}}
        },
        /pgfplots/box plot
    },
    box plot top whisker/.style={
        /pgfplots/error bars/draw error bar/.code 2 args={%
            \pgfkeysgetvalue{/pgfplots/error bars/error mark}%
            {\pgfplotserrorbarsmark}%
            \pgfkeysgetvalue{/pgfplots/error bars/error mark options}%
            {\pgfplotserrorbarsmarkopts}%
            \path ##1 -- ##2;
        },
        /pgfplots/table/.cd,
        y index=\pgfkeysvalueof{/pgfplots/box plot whisker top index},
        y error expr={
            \thisrowno{\pgfkeysvalueof{/pgfplots/box plot box top index}}
            - \thisrowno{\pgfkeysvalueof{/pgfplots/box plot whisker top index}}
        },
        /pgfplots/box plot
    },
    box plot bottom whisker/.style={
        /pgfplots/error bars/draw error bar/.code 2 args={%
            \pgfkeysgetvalue{/pgfplots/error bars/error mark}%
            {\pgfplotserrorbarsmark}%
            \pgfkeysgetvalue{/pgfplots/error bars/error mark options}%
            {\pgfplotserrorbarsmarkopts}%
            \path ##1 -- ##2;
        },
        /pgfplots/table/.cd,
        y index=\pgfkeysvalueof{/pgfplots/box plot whisker bottom index},
        y error expr={
            \thisrowno{\pgfkeysvalueof{/pgfplots/box plot box bottom index}}
            - \thisrowno{\pgfkeysvalueof{/pgfplots/box plot whisker bottom index}}
        },
        /pgfplots/box plot
    },
    box plot median/.style={
        /pgfplots/box plot,
        /pgfplots/table/y index=\pgfkeysvalueof{/pgfplots/box plot median index}
    },
    box plot width/.initial=1em,
    box plot x index/.initial=0,
    box plot median index/.initial=1,
    box plot box top index/.initial=2,
    box plot box bottom index/.initial=3,
    box plot whisker top index/.initial=4,
    box plot whisker bottom index/.initial=5,
}

\usepackage{framed}
\usepackage{booktabs}
\usepackage{caption}
\usepackage{float}
\usepackage{titlesec}
\usepackage{capt-of}
\definecolor{darkgreen}{rgb}{0.1, 0.5, 0.2}
\usepackage{array}
\usepackage{arydshln}
\newtheorem{definition}{Definition}

\definecolor{gray2}{rgb}{0.6,0.6,0.6}
\definecolor{lightgray2}{rgb}{0.8,0.8,0.8}
\title{An experimental approach: Converting verbal expressions to numerical scales} 
\author{Zsombor Sz\'adoczki$^{1,2,*}$, S\'andor Boz\'oki$^{1,2}$, L\'aszl\'o Sipos$^{3,4}$, Zs\'ofia Galambosi$^3$}
\date{}
\begin{document}
\pagenumbering{arabic}

\maketitle
\begin{center}
$*$ Corresponding author, 1111 Kende u. 13-17., Budapest, Hungary;\\ Email: szadoczki.zsombor@sztaki.hu\\
\bigskip
$^{1}$ Research Group of Operations Research and Decision Systems, \\
Research Laboratory on Engineering \& Management Intelligence \\
HUN-REN Institute for Computer Science and Control (SZTAKI), 1111 Kende u. 13-17., Budapest, Hungary;\\ Email: szadoczki.zsombor@sztaki.hu, bozoki.sandor@sztaki.hu\\
\bigskip
$^{2}$ Department of Operations Research and Actuarial Sciences \\
Corvinus University of Budapest, 1093 Fővám tér 8., Budapest, Hungary \\
\bigskip
$^3$ Hungarian University of Agriculture and Life Sciences,\\ Institute of Food Science and Technology, \\ Department of Postharvest, Commercial and Sensory Science, 1118, Villányi út 29-43., Budapest, Hungary;\\ Email: sipos.laszlo@uni-mate.hu, galambosi.zsofia@uni-mate.hu\\

\bigskip
$^4$ Institute of Economics,\\ HUN-REN Centre of Economic and Regional Studies, 1097 Tóth Kálmán utca 4., Budapest, Hungary;

\end{center}

\newpage
\renewcommand{\baselinestretch}{1.5} \normalsize

\begin{abstract}

\noindent
One of the key issues in decision problems is the selection and use of the appropriate response scale. In this paper verbal expressions are converted into numerical scales for a subjective problem instance. The main motivation for our research was that verbal values in decision tasks are often mechanically converted into numbers, which thus typically do not fully represent the respondent's true evaluation. In our experiment, we conducted a color selection test with 462 subjects by testing six colors (red, green, blue, magenta, turquoise, yellow) defined from the Color Namer database on color-calibrated tablets in ISO standardized sensory test booths of a sensory laboratory. The colors were evaluated both in a pairwise comparison matrix (indirect ranking with four-item verbal category scale) and on a direct scoring basis. We determined scales that provide the closest results on average and individually to the direct scoring, based on the eigenvector and the logarithmic least squares methods. All results show that the difference between verbal expressions is much smaller than the one used by most of the common numerical scales. The respondents' inconsistency was also analyzed, even with a repeated question regarding their preference between a given pair of colors. It is shown that most decision makers answer fairly similarly for the second time, but there can be significant (even ordinal) differences. The respondents whose answers are further from the original tend to be more inconsistent in general.


\end{abstract}

\noindent \textbf{Keywords}: Decision analysis, Pairwise comparisons, AHP, Scales, Empirical pairwise comparison matrix 

\renewcommand{\baselinestretch}{1.5} \normalsize

\section{Introduction}
\label{sec:1}

Pairwise comparisons are popular in decision making problems \citep{Triantaphyllou2000,Greco2025}, preference modelling \citep{DavidsonFarquhar1976}, as well as in sports \citep{Csato2021}. One of the most frequently used multicriteria decision making (MCDM) methods is the Analytic Hierarchy Process (AHP) proposed by Saaty \citep{Saaty1977,Saaty}, which applies pairwise comparison matrices (PCMs). The appropriate entry of a (multiplicative) PCM shows the decision maker's relative preference of the alternative in the row over the alternative in the column. 

Empirical pairwise comparison matrices are collected from real or experimental decision problems \citep{Bozoki2013,CavalloIshizaka2023}. They can substantially differ from simulated or specifically chosen example matrices that are mainly considered in the existing literature \citep{Szadoczki2022}. Thus, their analysis can lead to important consequences and practical recommendations in decision support.

In the classical AHP, the so-called Saaty-scale or fundamental scale---i.e., the discrete scale of $\{1/9,1/8, \ldots,$ $1/2,1,2,\ldots,8,9 \}$---is applied to make each comparison between the pairs of alternatives. However, these values are converted from the verbal assessments of the decision makers chosen from an appropriate discrete list of possibilities.

The data provided by empirical pairwise comparison matrices and alternative methods of preference elicitation, such as direct scoring, can be used to verify in what extent does a given scale corresponds to the decision makers' judgments.

In this paper, we conduct sensory-based experiments with pairwise comparisons to address the discrepancy between the commonly used scales and reality. The preferences of university students regarding different colors are collected via pairwise comparisons based on a four-item verbal category scale as well as direct scoring. The empirical distributions of the comparisons calculated from direct evaluation scores, and the verbal category scale applied for pairwise comparisons are analyzed to determine the optimal numerical scale that should be used between the conversion of verbal and numerical scales of pairwise comparisons in the given problem. Both the individually optimal scales, and the one that is the best on average are determined using the Eigenvector and the Logarithmic Least Squares Methods based on the Euclidean distances calculated with the direct evaluation scores. One of the questions regarding the preference between a pair of colors was repeated at the end of the experiment that also led to interesting findings regarding the inconsistency of the individuals.

The rest of the paper is organized as follows. Section \ref{sec:2} contains a literature review related to converting verbal scales to numerical ones in decision making problems. The methodology of the questionnaires and their evaluation are detailed in Section \ref{sec:3}. Section \ref{sec:4} presents the main results, the relation of direct scoring and pairwise comparisons, as well as the analysis of the respondents' inconsistency. Finally, Section \ref{sec:5} concludes and raises further research questions.

\section{Literature review}
\label{sec:2}

Empirical pairwise comparison matrices can significantly differ from simulated ones, and form an essential source of research. \cite{Bozoki2013} examine 454 PCMs (with 227 decision makers) and analyze the effect of the size, the question ordering, and the type of the problem (subjective or objective) on the inconsistency of the matrices. \cite{CavalloIshizaka2025} compare the precision of more than 100 respondents when they evaluate the area of nine geometric figures and five distances between cities applying pairwise comparisons, direct scoring, and the best-worst method \citep{Rezaei2015}. 

\cite{Huizingh1997} applies a laboratory study with 180 participants to compare the results of the Analytic Hierarchy Process using numerical values and verbal expressions converted by the fundamental scale ($\{1,2,\ldots,9\}$). The results show that the fundamental scale tends to overestimate the differences in preference, but the numerical mode only shows slightly better results (regarding inconsistency and the finally chosen alternative).

\cite{Poyhonen} perform a comparative study applying an experiment on areas of figures with different shapes. They compare three different scales (linear (fundamental), $9/9-9/1$, and balanced). The main finding of the experiment is that the perceived meaning of the verbal expressions varies from one decision maker to the other, and also depends on the set of elements involved in the comparison. They show that alternative
numerical scales yield more accurate estimates than the usual 1-to-9 scale, and these also reduce the inconsistency of the
comparison matrices measured by the consistency measure (CM) that is independent from the choice of the scale \citep{Salo1993,Salo1995}.

\cite{Elliott2010} examines three different scale types, and based on the answers of 64 respondents, finds that none of these scales can accurately capture the preferences of all individuals. On a real engineering decision problem it is shown that the choice of scale can also affect the preferred decision. \cite{Finan1999} apply experiments to argue for the calibration of the applied scales, particularly proposing a geometric scale based on a single parameter.

\cite{Dong2013} propose to generate the numerical scale of each individual based on a linguistic model in order to handle the problem caused by the fact that the individuals understand the same verbal scale in a different manner. They apply nonlinear programming to determine the numerical values based on five different scales known in the literature.

\cite{Rokou2014} show that the problem of choosing the appropriate scale is essential in group decision making. They present a group calibration process that adjusts individuals’ preferences based on their answers on a set of standardized questions. \cite{Meesariganda2017} compare eight different numerical scales. First, they apply the verbal scales to compare alternatives with known measures, e.g., the surface of figures. Then, based on the Euclidean distances, the best matching scale representing the real values is selected for each individual. \cite{Ahmed2023} propose a scale individualization approach based on compatibility. \cite{Dutta2025} recommend an approach comprising two steps to identify the most suitable personalized numerical scale for individuals.

\cite{Koczkodaj2017} propose a method to reduce the number of rating scale items without the loss of predictability, while \cite{Koczkodaj2019} apply differential evolution optimization to improve rating scale predictability.

\cite{Franek2014} compare eight different scales regarding the resulting priorities and consistency on an example problem. Results suggest that judgment
scales have a significant impact on criteria priorities but not on the ranking of criteria. The inconsistency varies among applied judgment scales for which the appropriate random indices have been calculated by the authors (see \citet[Table 7]{Franek2014}) to determine the suitable $CR$ inconsistency values. 

\cite{Goepel2019} proposes three different scales (generalized balanced, adaptive, and adaptive-balanced scales) also depending on the number of examined criteria, and provides a thorough comparative analysis as well.

\cite{CavalloIshizaka2023} conduct experiments, where the decision makers asked to compare the distance of different cities. Three problems are considered with 52 PCMs each, where eight different scales examined earlier in the literature are compared to each other. They find that the inverse linear scale provides the best results on average in every problem.

Table ~\ref{tab:1} contains some of the most popular scales proposed in the literature to assess verbal statements for pairwise comparison matrices. It is worth mentioning that even though all of these scales are presented in a general framework, they are used with a given parameter combination (shown in the ‘Usual choice of parameters' column in Table  ~\ref{tab:1}) in the vast majority of articles, and in practice as well.

Most of the studies comparing different scales focus on objective problems, e.g., evaluating the size of figures. They also deal with the usual parameter combinations of the different scales. This paper examines a large number of scales using a subjective problem, i.e., comparing different colors, which could be significantly different from objective tasks.

\begin{table}[H] \centering
\scriptsize
\renewcommand*{\arraystretch}{2}
\begin{tabular}{@{}ccccc>{\centering\arraybackslash}p{2.5cm}@{}}\toprule
	\multirow{ 2}{*}{\shortstack[c]{Name}} & \multirow{ 2}{*}{\shortstack[c]{Description}} & \multirow{ 2}{*}{\shortstack[c]{Parameters}} & \multirow{ 2}{*}{\shortstack[c]{Usual choice\\of parameters}} & \multirow{ 2}{*}{\shortstack[c]{Scale values}} & \multirow{ 2}{*}{\shortstack[c]{Reference}} \\
    
    &&&&&\\\hline
    
	Linear & $\alpha\cdot x$& $\alpha>0$, \ $x\in\{1,2,\ldots,9\}$ & $\alpha=1$ & $\{1,2,\ldots,9\}$ &\cite{Saaty1977} \\
    
    Affine & $\alpha\cdot x + \beta$& \multirow{ 2}{*}{\shortstack[c]{ $\alpha>0, \ \beta>0$,\\$x\in\{1,2,\ldots,9\}$}} & \multirow{ 2}{*}{\shortstack[c]{$\alpha=0.1$, \\ $\beta=1$}} & $\{1.1,1.2,\ldots,1.9\}$ &\cite{Saaty1987} \\
    
     & &  & &  \\
     
	Power & $x^{\alpha}$ &$\alpha>1$, \ $x\in\{1,2,\ldots,9\}$ &$\alpha=2$ & $\{1,4,9,\ldots,81\}$&\cite{HarkerVargas1987} \\
    
    Root & $\sqrt[\alpha]{x}$ &$\alpha>1$, \  $x\in\{1,2,\ldots,9\}$ & $\alpha=2$ & $\{1,\sqrt{2},\sqrt{3},\ldots,3\}$ &\cite{HarkerVargas1987} \\
    
    Geometric & $\alpha^{x-1}$ &\multirow{ 2}{*}{\shortstack[c]{$\alpha>1$, \  $x\in\{1,2,\ldots,9\}$,\\$x\in\{1,1.5,\ldots,4\}$, etc.}} & \multirow{ 2}{*}{\shortstack[c]{$\alpha=\sqrt{2}$ \\ or $\alpha=2$}} & \multirow{2}{*}{\shortstack[c]{$\{1,\sqrt{2},2,\ldots,16\}$\\or $\{1,2,4,\ldots,256\}$}} &\cite{Lootsma1989} \\
    
    &&& \\
    
    \multirow{ 2}{*}{\shortstack[c]{Inverse linear\\($9/9-9/1$)}}& $\frac{9}{10-x}$ & $x\in\{1,2,\ldots,9\}$ & - & $\{1,1.13,1.29\ldots,9\}$ & \cite{MaZheng1991} \\
    
   
    Asymptotic& $e^{tanh^{-1}\left(\frac{\sqrt{3}(x-1)}{14}\right)}$ & $x\in\{1,2,\ldots,9\}$ & - & $\{1,1.13,1.29\ldots,13.93\}$&\cite{Donegan1992} \\
    
    Balanced& $\frac{x}{1-x}\left( = \frac{9+y}{11-y}\right)$ & \multirow{ 2}{*}{\shortstack[c]{$x\in\{0.5,0.55,\ldots,0.9\}$,\\$y\in\{1,2,\ldots,9\}$}} & - & $\{1,1.22,1.5\ldots,9\}$ & \cite{Salo1997} \\
    
    Balanced power& $9^{\frac{x-1}{n-1}}$ & $x\in\{1,2,\ldots,n\}$ & $n=9$ & $\{1,1.32,1.73\ldots,9\}$ & \cite{Elliott2010} \\
    
    Logarithmic& $log_{\alpha}\left(x+\alpha-1\right)$ &$\alpha>1$, \ $x\in\{1,2,\ldots,n\}$ & $\alpha=2$ & $\{1,1.58,2\ldots,3.32\}$ &\cite{Ishizaka2011} \\

    Koczkodaj& $1+\frac{x-1}{n-1}$ & $x\in\{1,2,\ldots,n\}$ & $n=9$ & $\{1,1.125,1.25\ldots,2\}$ &\cite{Koczkodaj2016} \\
    
	\bottomrule
	\end{tabular}
 \caption{Popular scales proposed in connection with pairwise comparison matrices.}
 \label{tab:1}
\end{table}

\section{Methodology}
\label{sec:3}

\subsection{Pairwise comparison matrices}
\label{sec:3.1}
In this study, the concept of pairwise comparison matrices and related tools are widely used. Let us denote the number of alternatives by $n$. 
\begin{definition}[Pairwise comparison matrix (PCM)]
 The $n\times n$ matrix $\mathbf{A}=[a_{ij}]$ is a pairwise comparison matrix, if it is
 \begin{itemize}
     \item positive ($a_{ij}>0$ $\forall$  $i$,  $j$) and
     \item reciprocal ($1/a_{ij}  = a_{ji}$ $\forall$  $i$,  $j$).
 \end{itemize} 
\end{definition}

One of the most commonly analyzed questions is the (in)consistency of PCMs.

\begin{definition}[Consistent PCM]
A PCM is consistent if  $a_{ik}=a_{ij}a_{jk} \hspace{0.2cm} \forall i,j,k$. If a PCM is not consistent, then it is called inconsistent.
\end{definition}

The level of inconsistency can be measured several ways \citep{Brunelli2018}, and there are different systems to validate such inconsistency metrics \citep{Brunelli2024}. However, in practice, still the Consistency Ratio \citep{Saaty1977} is the most commonly used index for these purposes. It is also related to the $10\%$ rule (or $0.1$ rule), which states that a decision maker should be asked to reconsider a PCM in the case when the corresponding $CR$ is larger than $10\%$.

\begin{definition}[Consistency Ratio ($\boldsymbol{CR}$)]
The consistency ratio ($CR)$ of an $n\times n$ PCM $\mathbf{A}$ is determined by Equation~\ref{eq:1}.
\begin{equation}
\label{eq:1}
    CR=\frac{CI}{RI},
\end{equation}
where $CI$ denotes the Consistency Index, which is defined by Equation~\ref{eq:2}.
\begin{equation}
\label{eq:2}
    CI=\frac{\lambda_{\max}-n}{n-1},
\end{equation}
where $\lambda_{\max}$ is the principal eigenvalue of matrix $\mathbf{A}$, and $RI$ is the Random Index, which is the average $CI$ obtained from a sufficiently large set of randomly generated PCMs of size $n$ using the fundamental scale ($\{1,2,\ldots,9\}$) and discrete uniform distribution.
\end{definition}

It is common to calculate a weight vector from a given PCM that represents the estimation of the underlying preferences of the decision maker. Two of the most commonly used weight calculation techniques are the Eigenvector Method \citep{Saaty1977}, and the Logarithmic Least Squares (Geometric Mean) Method \citep{Crawford1985}.

\begin{definition}[Eigenvector Method (EM)]
The weight vector $\mathbf{w}$ of an $n\times n$ PCM $\mathbf{A}$ determined by the EM is defined by Equation~\ref{eq:3}.
\begin{equation}
\label{eq:3}
    \mathbf{A} \mathbf{w}=\lambda_{\max}\mathbf{w},
\end{equation}
where the componentwise positive principal eigenvector $\mathbf{w}$ is unique up to a scalar multiplication.
\end{definition}

\begin{definition}[Logarithmic Least Squares Method (LLSM)]
The weight vector $\mathbf{w}$ of an $n\times n$ PCM $\mathbf{A}$ determined by the LLSM is given by Equation~\ref{eq:4}.
\begin{equation}
\label{eq:4}
\min_{\mathbf{w}}    \sum_{i=1}^n\sum_{j=1}^n \left(\ln(a_{ij})-\ln\left(\frac{w_i}{w_j}\right)\right)^2 ,
\end{equation}
where $w_i$ is the $i$th coordinate of $\mathbf{w}$.
\end{definition}

\subsection{Framework of the ratio scale generation and comparison}
\label{sec:3.2}

In an experiment, it is interesting to examine the scale that provides the best results on average. Where ‘the best' in our methodology means the closest (based on the Euclidean distance) to the objective values for an objective analysis, and the closest to the values calculated from other preference elicitation methods, such as the direct scoring, in the case of subjective problems. However, it also makes sense to find the optimal scale for each individual as well.

In order to answer these questions in our experiment, the EM and LLSM weight vectors of the individual pairwise comparison matrices have been calculated for several scales, applying the following parametrization of the used four-item verbal scale.

\begin{itemize}
    \item $1$ -- equally like;
    \item $S$ -- like it a little bit more;
    \item $M$ -- like it moderately more;
    \item $L$ -- like it much more.
\end{itemize}

Where $1<S<M<L$, thus a strict relation is assumed between the different levels, and the category representing indifference (‘equally like') must be converted to $1$ as a numerical value.

All parameter combinations were tested on a grid with a step size of $0.1$ in each variable, where the above assumptions hold, and $S\leq5$, $M\leq10$, and $L\leq 15$. That is a total of 236 880 parameter combinations.

\subsection{Conducted experiments}
\label{sec:3.3}
In the conducted color-choice test, the colors selected were those that were clearly distinguishable on the subjects' mental maps: red, green, blue, magenta, turquoise and yellow. The RGB coordinates of these color names can be found in the Color Namer database (\url{https://colornaming.net/#colour-namer}): red (R:189, G:62, B:57), green (R:90, G:151, B:90), blue (R:84, G:110, B: 183), magenta (R: 179, G: 55, B: 151), turquoise (R: 63, G: 185, B: 177), and yellow (R: 227, G: 203, B: 78) \citep{MylonasMacDonald2017}. The tests were conducted on color-calibrated tablets (Samsung Galaxy Tab A 2018, 10.5) with identical test geometry and maximum brightness settings in the sensory test booths of the Sensory Laboratory, Institute of Food Science and Technology, Hungarian University of Agricultural and Life Sciences. The laboratory conforms to the International Standard for the design requirements for sensory test  \citep{ISO8589}.

In the test, each color was initially compared with every other color in pairs based on their stated preferences, resulting in a total of 15 pairs (calculated as $n(n-1)/2$). In the first phase, participants were first asked to indicate their preferred color in pairs based on their preference, and then to rate their level of preference using a four-item verbal category scale: equally like, like it a little bit more, like it moderately more, like it much more (see Figure~\ref{fig:Appendix} in Appendix~\ref{append:A}). In the second phase of the experiment, the participants were asked to rate the colors on an eleven-point scale, with 0 indicating a strong dislike and 10 indicating a strong preference (see Figures~\ref{fig:Appendix2} and \ref{fig:Appendix3} in Appendix~\ref{append:A}). Consequently, the colors were evaluated both on a pairwise comparison (indirect ranking) and on an individual basis (direct ranking). In the third phase, participants were asked a series of socio-demographic questions. These included queries regarding gender, age, and county. To assess the consistency of the responses, the question regarding the preference of the secondly asked pair of colors was repeated at the conclusion of the test in reverse order.




\section{Results}
\label{sec:4}

There were a total of 462 respondents in our color-choice test. However, before the calculations, it was necessary to clean the dataset as follows.

\begin{enumerate}
    \item The results of the individuals who did not use the whole verbal scale (e.g., there were no pair of colors for which they liked one a little bit more than the other) were removed. This could have caused problems when the individually best scales were determined as there would have been a number of optimal cases. For instance, in the previous example $S$ is not used, thus, $M$ and $L$ are optimized, and $S$ can take any value that satisfies the assumptions. This results in the removal of the answers of 98 respondents.
    \item It is easy to compare the numerical scale values with the ratio of the direct scoring results. However, it was possible to determine 0 as the score of a color. The answers of the 63 respondents who used this value were also removed from the calculations.
\end{enumerate}

After that a total of 301 PCMs and direct scoring results were evaluated. It is also worth mentioning that we did calculate the results without excluding these respondents (using transformations for the second point), and all the main findings were the same.

\subsection{Optimal scales}
\label{sec:4.1}

\begin{figure}[H]
    \centering

\begin{tikzpicture}[scale=0.53]

\begin{axis} [ybar, height=10cm,width=21cm,max space between ticks=75pt,
        try min ticks=5,bar width=12pt,xtick={0,1,2,3,4,5,6,7,8,9,10},ymin=0,xmin=0,xlabel={Ratio of the direct evaluation scores},ylabel={Frequency}]{
ymin=0,
ymax=200,
xmin=0,
xtick={0,1,2,3,4,5,6,7,8,9,10},
yticklabel style={
        /pgf/number format/fixed
},
xticklabel style={
        /pgf/number format/fixed,
},
}
\addplot coordinates {
    (0.25,4)
(0.5,26)
(0.75,50)
(1,177)
(1.25,65)
(1.5,61)
(1.75,9)
(2,13)
(2.25,0)
(2.5,7)
(2.75,1)
(3,2)
(3.25,0)
(3.5,2)
(3.75,0)
(4,2)
(4.25,0)
(4.5,0)
(4.75,0)
(5,1)
(5.25,0)
(5.5,0)
(5.75,0)
(6,1)
(6.25,0)
(6.5,0)
(6.75,0)
(7,0)
(7.25,0)
(7.5,0)
(7.75,0)
(8,0)
(8.25,0)
(8.5,0)
(8.75,0)
(9,0)
(9.25,0)
(9.5,0)
(9.75,0)
(10,0)
};
\end{axis}

\begin{axis} [yshift=-10cm, ybar, height=10cm,width=21cm,max space between ticks=75pt,
        try min ticks=5,bar width=12pt,xtick={0,1,2,3,4,5,6,7,8,9,10},ymin=0,xmin=0,xlabel={Ratio of the direct evaluation scores},ylabel={Frequency}]{
ymin=0,
ymax=200,
xmin=0,
xtick={0,1,2,3,4,5,6,7,8,9,10},
yticklabel style={
        /pgf/number format/fixed
},
xticklabel style={
        /pgf/number format/fixed,
},
}
\addplot coordinates {
(0.25,0)
(0.5,17)
(0.75,61)
(1,308)
(1.25,302)
(1.5,251)
(1.75,128)
(2,120)
(2.25,5)
(2.5,57)
(2.75,20)
(3,37)
(3.25,0)
(3.5,17)
(3.75,0)
(4,16)
(4.25,0)
(4.5,2)
(4.75,0)
(5,16)
(5.25,0)
(5.5,0)
(5.75,0)
(6,6)
(6.25,0)
(6.5,0)
(6.75,0)
(7,5)
(7.25,0)
(7.5,0)
(7.75,0)
(8,2)
(8.25,0)
(8.5,0)
(8.75,0)
(9,1)
(9.25,0)
(9.5,0)
(9.75,0)
(10,5)

};
\end{axis}

\begin{axis} [yshift=-20cm, ybar, height=10cm,width=21cm,max space between ticks=75pt,
        try min ticks=5,bar width=12pt,xtick={0,1,2,3,4,5,6,7,8,9,10},ymin=0,xmin=0,xlabel={ Ratio of the direct evaluation scores},ylabel={Frequency}]{
ymin=0,
ymax=200,
xmin=0,
xtick={0,1,2,3,4,5,6,7,8,9,10},
yticklabel style={
        /pgf/number format/fixed
},
xticklabel style={
        /pgf/number format/fixed,
},
}
\addplot coordinates {
(0.25,3)
(0.5,13)
(0.75,25)
(1,185)
(1.25,230)
(1.5,288)
(1.75,158)
(2,169)
(2.25,33)
(2.5,79)
(2.75,41)
(3,52)
(3.25,0)
(3.5,35)
(3.75,0)
(4,22)
(4.25,0)
(4.5,6)
(4.75,0)
(5,27)
(5.25,0)
(5.5,0)
(5.75,0)
(6,7)
(6.25,0)
(6.5,0)
(6.75,0)
(7,7)
(7.25,0)
(7.5,0)
(7.75,0)
(8,8)
(8.25,0)
(8.5,0)
(8.75,0)
(9,3)
(9.25,0)
(9.5,0)
(9.75,0)
(10,3)

};
\end{axis}

\begin{axis} [yshift=-30cm, ybar, height=10cm,width=21cm,max space between ticks=75pt,
        try min ticks=5,bar width=12pt,xtick={0,1,2,3,4,5,6,7,8,9,10},ymin=0,xmin=0,xlabel={Ratio of the direct evaluation scores},ylabel={Frequency}]{
ymin=0,
ymax=200,
xmin=0,
xtick={0,1,2,3,4,5,6,7,8,9,10},
yticklabel style={
        /pgf/number format/fixed
},
xticklabel style={
        /pgf/number format/fixed,
},
}
\addplot coordinates {
(0.25,0)
(0.5,6)
(0.75,14)
(1,56)
(1.25,129)
(1.5,169)
(1.75,128)
(2,208)
(2.25,43)
(2.5,104)
(2.75,48)
(3,78)
(3.25,0)
(3.5,86)
(3.75,0)
(4,54)
(4.25,0)
(4.5,42)
(4.75,0)
(5,51)
(5.25,0)
(5.5,0)
(5.75,0)
(6,13)
(6.25,0)
(6.5,0)
(6.75,0)
(7,15)
(7.25,0)
(7.5,0)
(7.75,0)
(8,25)
(8.25,0)
(8.5,0)
(8.75,0)
(9,27)
(9.25,0)
(9.5,0)
(9.75,0)
(10,28)

};
\end{axis}

\node [fill=white,draw=white] at (10,6){ $\boldsymbol{1\sim}$ \textbf{Equally like}};

\node [fill=white,draw=white] at (10,-4){ $\boldsymbol{S\sim}$ \textbf{Like it a}};

\node [fill=white,draw=white] at (10,-5){\textbf{little bit more}};

\node [fill=white,draw=white] at (10,-14){ $\boldsymbol{M\sim}$ \textbf{Like it}};

\node [fill=white,draw=white] at (10,-15){\textbf{moderately more}};

\node [fill=white,draw=white] at (10,-24){ $\boldsymbol{L\sim}$ \textbf{Like it}};

\node [fill=white,draw=white] at (10,-25){\textbf{much more}};
 
\end{tikzpicture}

 \caption{ \footnotesize The histograms of the $1$, $S$, $M$, and $L$ values from the direct evaluation, respectively, i.e., the histograms of the ratios of directly evaluated values, where the pairwise comparison resulted as ‘equally like' (indifference), ‘like it a little bit more', ‘like it moderately more', and ‘like it much more' for the given two colors, respectively.}
    \label{fig:1SMLhist}

    \end{figure}
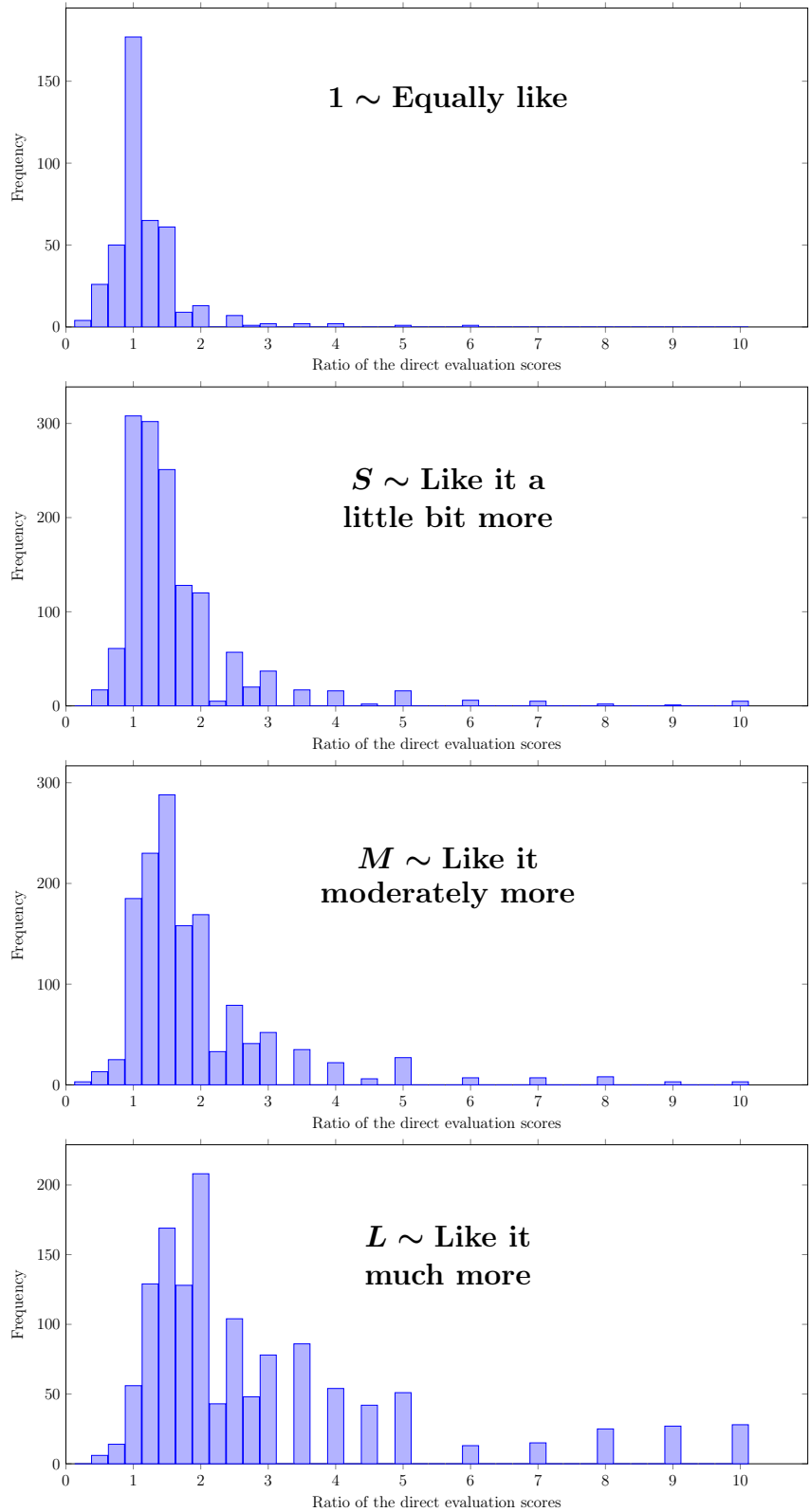

Figure~\ref{fig:1SMLhist} shows the histograms of the $1$, $S$, $M$, and $L$ parameters from the direct evaluations, respectively, i.e., the histograms of the ratios of directly evaluated values, where the pairwise comparison resulted as ‘equally like’ (indifference), ‘like it a little bit more’, ‘like it moderately more’, and ‘like it much more’ for the given two colors, respectively.

One can see that in the first histogram which is representing the indifference between the two colors, $1$ is the most common value of the ratio of direct evaluation scores. The distribution is more or less symmetric to that in the sense that the ratio value is smaller than $1$ approximately the same number of times as it is greater than $1$.

Then, as we go on to the more extreme preference values (like it a little bit more, moderately more, and much more) the common values for the ratio of direct evaluation scores starting to be larger and larger as one would expect. However, interestingly enough, there are still smaller than $1$ ratios for the most extreme case as well, showing the inconsistency of the respondents.

\begin{figure}[H]
    \centering

\begin{tikzpicture}
\begin{axis}[
    xlabel=$S$,
    ylabel=$M$,
    zlabel={$L$},
    grid=major,
    xmin =1,
    xmax = 4,
    ymin = 1,
    ymax = 4,
    zmin = 1,
    zmax = 4,
    xtick= {1,2,3,4},
    ytick= {2,3,4},
    ztick= {2,3,4},
    title=The optimal scales of individuals for the EM,
            colorbar,
            colormap/hot
]

\addplot3[scatter, only marks, visualization depends on = {\thisrowno{3}*0.2+0.2 \as \perpointmarksize},
    scatter/@pre marker code/.append style={/tikz/mark size=\perpointmarksize}] table [x index=0, y index=1, z index=2,scatter src=\thisrowno{3}, col sep=space] {Figures/Final_EM_Individual444.dat};

 \addplot3 [color=black, opacity=0.5] coordinates {(4.45,0,0) (4.45,0,4.6)};

  \addplot3 [color=black, opacity=0.5] coordinates {(0,-1,5.3) (4.45,0,4.6)};

  \addplot3 [color=black, opacity=0.5] coordinates {(4.45,0,4.6) (5,4.35,5.35)};

\end{axis}
\end{tikzpicture}
\caption{The optimal scales of individuals based on the Eigenvector Method, $1.1\leq S\leq3.8$, \ $1.2\leq M\leq3.9$, \ $1.3\leq L\leq4$, \ and $S < M < L$.}
\label{fig:EMIndividual}
\end{figure}
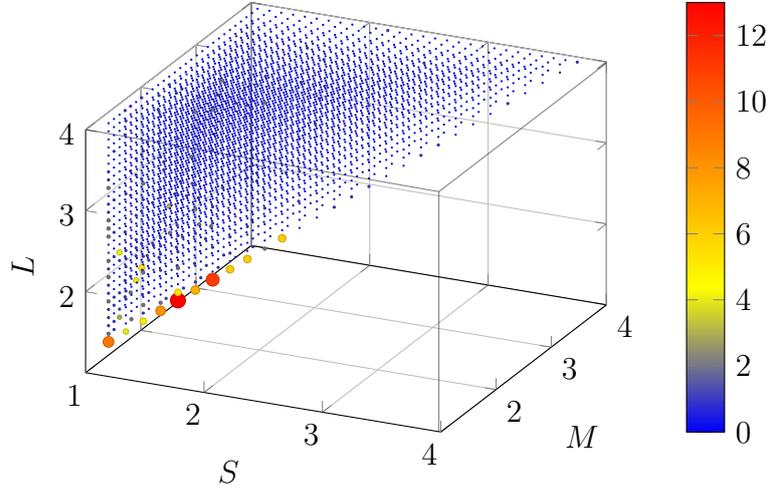

\begin{figure}[H]
\centering
    \includegraphics[width=0.55\paperwidth]{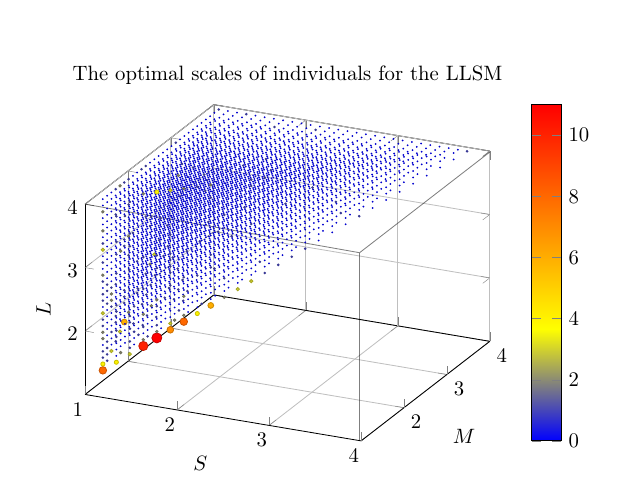}
 \caption{The optimal scales of individuals based on the Logarithmic Least Squares Method, $1.1\leq S\leq3.8$, \ $1.2\leq M\leq3.9$, \ $1.3\leq L\leq4$, \ and $S < M < L$.}
\label{fig:LLSMIndividual}
    \end{figure}

Here we only present the results for those parameter combinations, where $S,M,L\leq4$, while all the other assumptions hold, as it turned out that the scales including larger values are almost never optimal for any individuals.
Figures~\ref{fig:EMIndividual} and \ref{fig:LLSMIndividual} show the number of individuals for whom a given scale has been optimal based on the Euclidean distance using the Eigenvector and the Logarithmic Least Squares weight calculation techniques, respectively. The greater number of individuals is highlighted by the color and the size of the given point as well.

As one can see, the results of the two weight calculation methods are similar, and the most important conclusions are the same for both. All the commonly optimal scales include smaller than $3$ parameter values for even $L$. The scale that results in the smallest average Euclidean distance uses the parameter combination of $S=1.5$, $M=1.7$, and $L=2$ for both weight calculation techniques.

However, it is important to highlight the fact that the parameter combinations that are relatively close to the one that is proved to be optimal on average according to the Euclidean distance, all provide close-to-optimal values. Meanwhile, the average distance for the scales that include larger values ($\geq3$) are far from being optimal.

\subsection{The inconsistency of participants}
\label{sec:4.2}

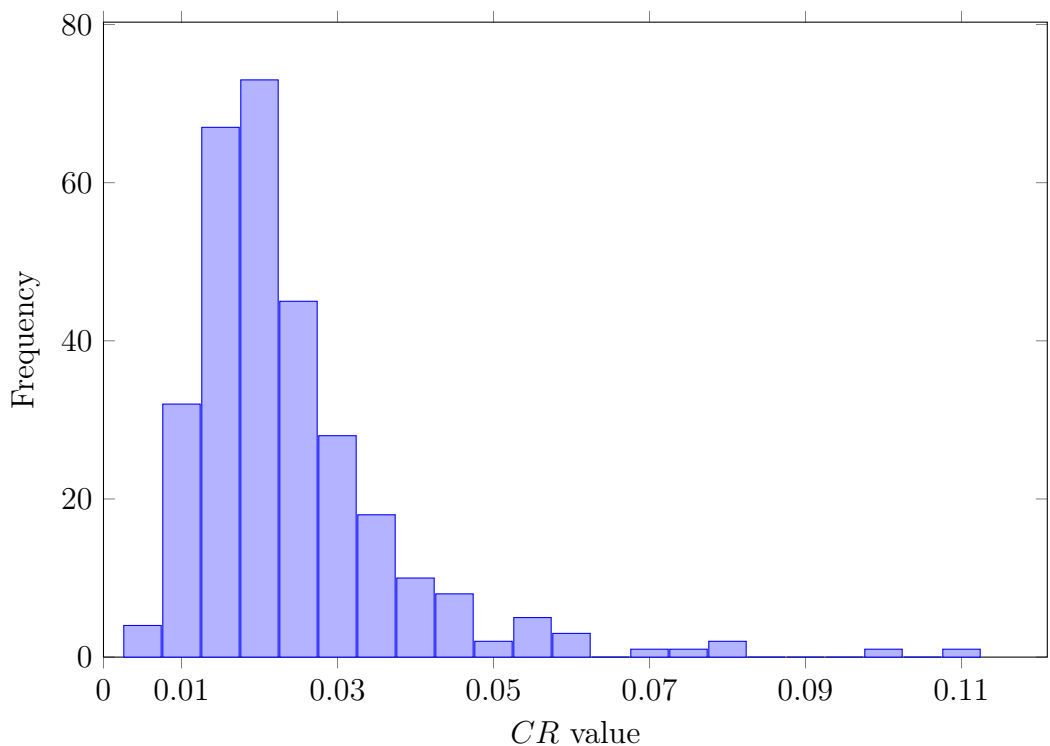
\begin{figure}[H]
    \centering

\begin{tikzpicture}
 
\begin{axis} [ybar, height=10cm,width=14cm,max space between ticks=75pt,
        try min ticks=5,bar width=14pt,xtick={0,0.01,0.03,0.05,0.07,0.09,0.11
},ymin=0,xmin=0,xlabel={$CR$ value},ylabel={Frequency},xticklabel style={
        /pgf/number format/fixed}]{
ymin=0,
ymax=100,
xmin=0,
xtick={0,0.01,0.03,0.05,0.07,0.09,0.11
},
yticklabel style={
        /pgf/number format/fixed
},
xticklabel style={
        /pgf/number format/fixed,
},
}
\addplot coordinates {
(0.005,4)
(0.01,32)
(0.015,67)
(0.02,73)
(0.025,45)
(0.03,28)
(0.035,18)
(0.04,10)
(0.045,8)
(0.05,2)
(0.055,5)
(0.06,3)
(0.065,0)
(0.07,1)
(0.075,1)
(0.08,2)
(0.085,0)
(0.09,0)
(0.095,0)
(0.1,1)
(0.105,0)
(0.11,1)
};
\end{axis}
 
\end{tikzpicture}

 \caption{The histogram of the $CR$ values based on the $S=1.5$, $M=1.7$, and $L=2$ scale.}
    \label{fig:InconsistencyFinal}

    \end{figure}

 One can see the histogram of the respondents' inconsistencies measured by the $CR$ index in Figure~\ref{fig:InconsistencyFinal}.

 However, the $CR$ consistency index heavily builds on the fundamental scale with the generation of the $RI$. For a $6\times6$ matrix, the original $RI$ is 1.249 calculated by simulations in several studies (see, for instance, \cite{Bozoki2008} or \cite{Agoston2022}). We ran a simulation with 10 million randomly generated PCMs using the scale which resulted in the smallest average Euclidean distance ($S=1.5$, $M=1.7$, and $L=2$). The calculated modified random index (let us denote it by $\widehat{RI}$) is 0.09224. Thus, a modified $\widehat{CR}$ can be calculated using $\widehat{RI}$, where the original $CR$ index is multiplied by $\sim 13.5$.

This shows the importance of the numerical scale in evaluating inconsistency as well. However, we do not recommend the application or generalization of the common $10\%$ rule in this case. The problem is highly subjective, and of a larger size ($6\times 6$ PCM). Furthermore, the number of options is also much smaller than the ones applied in the fundamental scale resulting in differences of the measured inconsistency.

 To assess the consistency of the respondents in a different way, the preference elicitation of the secondly asked pair was repeated at the conclusion of the test in reverse order.
 
\begin{figure}[H]
    \centering

\begin{tikzpicture}
	\pgfplotstableread[col sep=comma]{data.csv}\csvdata
	\pgfplotstabletranspose\datatransposed{\csvdata} 
	\begin{axis}[width=15cm,
		boxplot/draw direction = y,
		x axis line style = {opacity=1},
		axis x line = bottom,
            xmin = 0,
            xmax = 8,
		axis y line = left,
		enlarge y limits,
		ymajorgrids,
		xtick = {1, 2, 3, 4, 5, 6, 7},
		xticklabel style = {align=center, font=\small, below=2.5mm}, 
		xticklabels = {Same, $1$ step, $2$ steps, $3$ steps, $4$ steps, $5$ steps, $6$ steps},
		xtick style = {draw=none}, 
		ylabel = {$CR$},
            yticklabel style={
        /pgf/number format/fixed,
        /pgf/number format/precision=2
},
scatter/classes={
    d={mark=otimes*,draw=brown, fill=brown}},
	]
		\foreach \n in {1,...,6} {
			\addplot+[boxplot, fill, draw=black] table[y index=\n] {\datatransposed};
		}

        \addplot[scatter,only marks,
    scatter src=explicit symbolic]table[meta=label] {
x y label
7	0.019134459 d
    };
	\end{axis}

\node at (-1,-0.4){\footnotesize Distance from};
\node at (-1,-0.8){\footnotesize the original element};

\node at (-1,-1.4){\footnotesize Number of};
\node at (-1,-1.8){\footnotesize individuals};

\node at (1.65,-1.6){\footnotesize 124};
\node at (3.35,-1.6){\footnotesize 108};
\node at (5.05,-1.6){\footnotesize 39};
\node at (6.75,-1.6){\footnotesize 16};
\node at (8.45,-1.6){\footnotesize 8};
\node at (10.15,-1.6){\footnotesize 5};
\node at (11.85,-1.6){\footnotesize 1};
    
\end{tikzpicture}

 \caption{The $CR$ inconsistency of the individuals depending on the distance (measured in the steps of the verbal scale) of the repeated question's result from the original answer via Box plots (highlighting the minimum, first, second, and third quartile, and the maximum). The number of individuals included in the given category is displayed for each case as well. Instead of a Box plot, only a point representing the corresponding $CR$ is shown for the category with only one individual.}
    \label{fig:repeated}

    \end{figure}
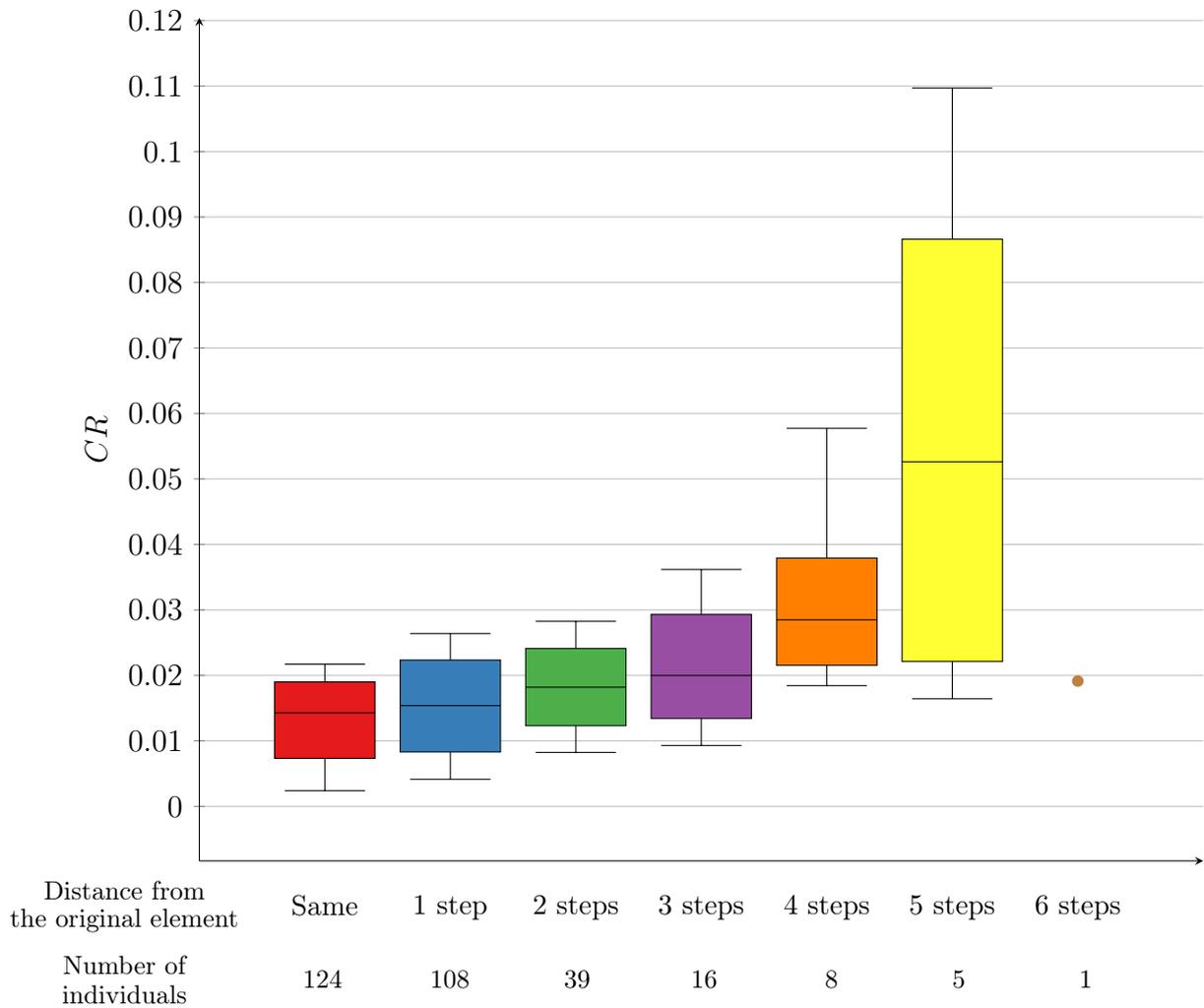

 Figure~\ref{fig:repeated} shows the $CR$ inconsistencies for the previously examined PCMs via Box plots for the possible categories of individuals depending on how far their repeated answer was from the original, as well as the number of individuals in the different categories.

   One can see that the measured inconsistency tends to be higher as the repeated answer is further from the original one. There are much less individuals who made totally different judgments for the second time compared to their original answer. More than 77\% of the respondents ($232/301$) made the same judgment for the second time as well, or only modified their answer with one step in the verbal scale.  However, we can even see one respondent, who reversed the most extreme preference between the pair. For that individual, only a point measuring the corresponding $CR$ is displayed on Figure~\ref{fig:repeated} instead of a Box plot.

  It is also worth mentioning that if there are at least $3$ steps in the four-item verbal scale between the original and the repeated preference, that means even the ordinal preference (which color is preferred) is changed at least to indifference. If the number of steps is even larger than $3$, then the preference must be reversed. However, the number of individuals whose decisions resulted in a difference of at least $3$ steps are a little bit less than 10\% of the respondents ($30/301$).

 There are two interesting findings based on the repeated preference elicitation.

 \begin{enumerate}
     \item Even though most of the respondents ($\sim 77\%$) answered almost the same way, there can be significant (even ordinal) differences in the results depending on when we ask the decision maker a given question.
     \item The individuals who are inconsistent regarding their (repeated) answer to a question tend to be more inconsistent in general.
 \end{enumerate}

\section{Conclusion and further research}
\label{sec:5}

In this paper, we analyzed the individually and on average optimal numerical scales related to a four-item verbal scale based on a color-choice test carried out on color-calibrated tablets in ISO standardized sensory laboratory and test booths \citep{ISO8589}. 462 respondents participated in the experiment providing verbal pairwise comparisons, as well as direct evaluation scores regarding the six colors (red, green, blue, magenta, turquoise, and yellow).

The Eigenvector and the Logarithmic Least Squares Methods were used to calculate the preference representing weight vectors from the pairwise comparison matrices. To compare the obtained weight vectors with the corresponding values calculated from the direct evaluation scores, the Euclidean distance was used. A total of 236 880 different numerical scales were considered.

It turned out that both on average and individually the numerical scales using smaller ($\leq 3$) numbers even for the extreme preference prevail in this subjective problem. The results do not depend on the used weight calculation technique. This is a significant finding, as it contradicts to both the common practice in Multicriteria Decision Making (especially the AHP), and to several ISO standards related to sensory tests, where verbal expressions are mechanically converted into great numbers (e.g., 9).

Following the results of transforming the verbal expressions of the category scale into numerical values, it is necessary to revisit some of the international standards, such as the ‘Guidelines for the use of quantitative response scales’ \citep{ISO4121}, ‘Selection and training of sensory assessors’ \citep{ISO8586}, ‘Guidelines for monitoring the performance of a quantitative sensory panel’ \citep{ISO11132}, and ‘General guidance for conducting hedonic tests with consumers in a controlled area’ \citep{ISO11136}.

We also analyzed the inconsistency of the respondents. It was emphasized that the scale itself heavily influences most of the inconsistency measures (for instance the $CR$ inconsistency index). The modified random index ($RI$) related to a modified $CR$ for the numerical scale that proved to be optimal was calculated. However, the generalization of the known $10\%$ rule is not recommended based on the different number of possibilities (four-item scale instead of nine-item one), the nature (subjective) and size ($6\times6$) of the problem.

A further inconsistency analysis was carried out in the paper, that included a repeated answer to a preference relation between a given pair of colors at the conclusion of the experiment. The ratio of respondents who provided answers that resulted in a sure ordinal difference in the preference (at least a distance of $3$ steps in the four-item verbal scale) was a little bit less than $10\%$. However, $\sim 77\%$ of the respondents answered almost the same way for the second time as well (at most a distance of $1$ step in the four-item verbal scale). Examining the inconsistency of respondents whose repeated answer was further from the original, it turned out that they tend to be more inconsistent in general.

Further research questions include analyzing the conversion between verbal and numerical scales on other subjective problems as well, e.g., applying further sensory tests (odour, flavor, etc.). A future research can also investigate the extent in which the used numerical scales can be personalized in several fields of decision making and in different practical problems as well.

The use of verbal category scales is widespread in preference studies, so the cognitive conversion of verbal categories into numbers is necessary.  Our specific results show how each verbal category can be converted to a numerical value in a given color-choice test task. However, in other tasks involving visual or other modality stimuli (taste, smell, texture), the numerical values of the same verbal categories need to be determined on a test-by-test, assessor-by-assessor basis, as this provides a reliable basis for mapping respondents' preferences more closely to reality. Furthermore, it is important to underline that there is currently very limited knowledge available on the stability of individual cognitive categories, including their dependence on the type of test item, the modality, the environment, the product being tested, etc. In future research, there will certainly be an increased role for the application and development of numerical conversion methods of verbal category scales.

\section*{Acknowledgements}
The project identified by EKOP-CORVINUS-24-4-080 was realized with the support of the National Research, Development, and Innovation Fund provided by the Ministry of Culture and Innovation, as part of the University Research Scholarship Program announced for the 2024/2025 academic year. The research was supported by the National Research, Development and Innovation Office under Grants FK 145838 and TKP2021-NKTA-01 NRDIO.

\bibliographystyle{apalike} 
\bibliography{main}
\addcontentsline{toc}{section}{References}

\newcounter{appendfigure}
\setcounter{appendfigure}{1}

\newcounter{append}
\renewcommand{\theappend}{A}

\section*{Appendix A}
\addcontentsline{toc}{section}{Appendix A}
\refstepcounter{append}
\label{append:A}

\begin{figure}[H]
\renewcommand{\thefigure}{A\arabic{appendfigure}}
\centering
    \includegraphics[width=0.5\paperwidth]{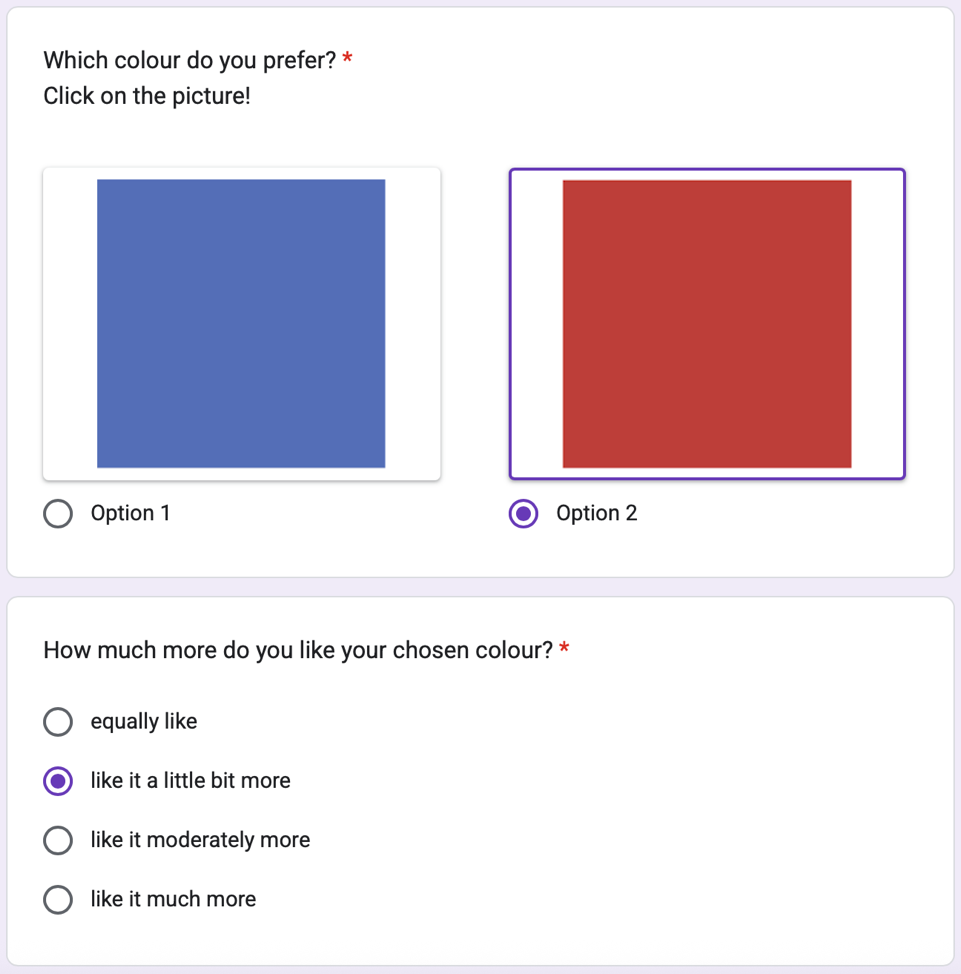}
 \caption{A representative example of the questions for pairwise comparisons used in our test.}
    \label{fig:Appendix}
    \setcounter{appendfigure}{2}
    \end{figure}

\begin{figure}[H]
\renewcommand{\thefigure}{A\arabic{appendfigure}}
\centering
    \includegraphics[width=0.5\paperwidth]{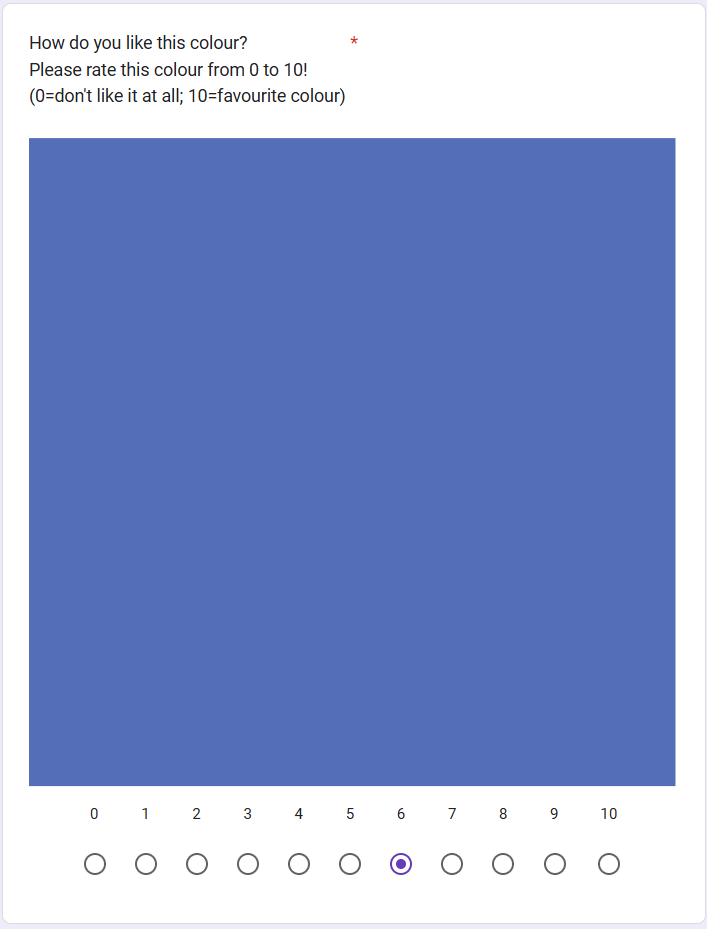}
 \caption{A representative example of the questions for direct evaluation used in our test.}
    \label{fig:Appendix2}
    \setcounter{appendfigure}{3}
    \end{figure}

\begin{figure}[H]
\renewcommand{\thefigure}{A\arabic{appendfigure}}
\centering
    \includegraphics[width=0.5\paperwidth]{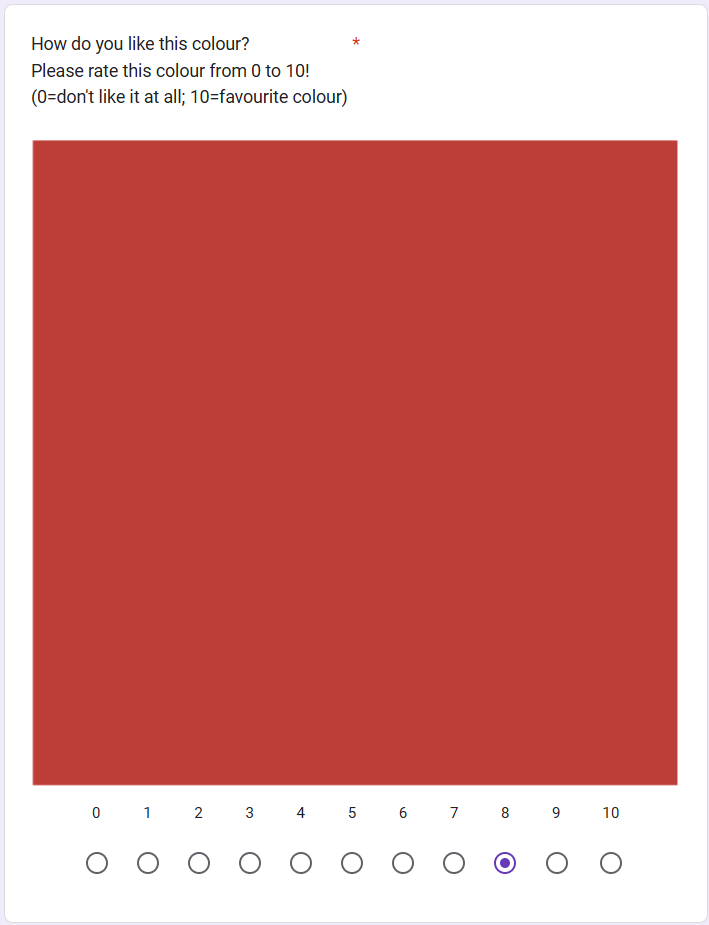}
 \caption{Another example of the questions for direct evaluation used in our test.}
    \label{fig:Appendix3}
    \setcounter{appendfigure}{4}
    \end{figure}

\end{document}